\documentclass[10pt]{article}    
\newtheorem{theorem}{Theorem}
\newtheorem{corollary}[theorem]{Corollary}
\newtheorem{remark}[theorem]{Remark}
\newtheorem{definition}[theorem]{Definition}

\usepackage{fancyhdr}
\title{\textbf{Local convergence analysis of the Gauss-Newton-Kurchatov method}}  

\author{Ioannis K. Argyros$^1$, Stepan Shakhno$^2$\\$^1$Department of Mathematics, Cameron University,\\ Lawton, USA, OK  73505; \\iargyros@cameron.edu,\\
\\$^2$Department of Theory of Optimal Processes, \\Ivan Franko National University of Lviv,\\ Lviv, Ukraine, 79000; \\stepan.shakhno@lnu.edu.ua} 
\begin{document}           
\maketitle                 

\begin{abstract}{\bf Abstract. }
We present a local convergence analysis of the Gauss-Newton-Kurchatov method  for solving nonlinear least squares problems with a decomposition of the operator.  The method uses the sum of the derivative of the differentiable part of the operator and the divided difference of the nondifferentiable part instead of computing  the full Jacobian. A theorem, which establishes the conditions of convergence, radius and the convergence order of the proposed method, is proved (Shakhno 2017). However, the radius of convergence is small in general limiting the choice of initial points. Using tighter estimates on the distances, under weaker hypotheses (Argyros et al. 2013), we provide an analysis of the Gauss-Newton-Kurchatov method with the following advantages over the corresponding results (Shakhno 2017): extended convergence region; finer error distances, and an at least as precise information on the location of the  solution. The numerical examples illustrate the theoretical results. 
\end{abstract}

\vspace*{0.1cm}

\vspace*{0.1cm}
\section{Introduction}
Let us consider the problem of finding an approximate solution of the nonlinear least squares problem  
\begin{equation} \label{EQ__1_} 
{\mathop{\min }\limits_{x\in {\bf {\rm R}}^{n} }} \frac{1}{2} F(x)^{{\bf \top }} F(x), 
\end{equation} 
where the residual function $F:D\subseteq {\bf {\rm R}}^{n} \to {\bf {\rm R}}^{m} $, $m\ge n$ is nonlinear in  $x$, $F$ is continously differentiable, and $D$ is an open convex set in ${\bf {\rm R}}^{n} $. 

A large number of problems in applied mathematics and also in engineering are solved by finding the solutions of problem (\ref{EQ__1_}). For example,  solving  overdetermined systems  of  nonlinear equations, estimating  parameters of physical processes  by measurement results, constructing nonlinear regressions models for solving  engineering,  problems dynamic systems, etc. The used solution methods are iterative -- when starting from one or several initial approximations a sequence is constructed that converges to a solution of the problems (\ref{EQ__1_}). 

Known methods of the Gauss-Newton type (Dennis et al. 1996; Ortega et al. 1970; Argyros 2008; Shakhno 2001) are used to solve the problem (\ref{EQ__1_}), which have derivatives of function in their iterative formulas. However, in practice, problems with calculations of derivative  arise. In this case, we can use iterative-difference methods (Argyros 2008; Ren et al. 2010, 2011; Shakhno et al. 1999, 2005) that do not require the calculation of the matrix of derivatives and often are not inferior  over the Gauss-Newton method at the order of convergence and the number of iterations. But sometimes the nonlinear function consists of differentiable and non-differentiable parts. Then a nonlinear least squares problem arises 
\begin{equation} \label{EQ__2_} 
{\mathop{\min }\limits_{x\in {\bf {\rm R}}^{n} }} \frac{1}{2} (F(x)+G(x))^{{\bf \top }} (F(x)+G(x)),      
\end{equation} 
where the residual function $F+G:D\subseteq {\bf {\rm R}}^{n} \to {\bf {\rm R}}^{m} $, $m\ge n$,  is nonlinear in $x$, $F$ is continously differentiable, $G$ is continous function, differentiability of which, in general, is not assumed,  and $D$ is an open convex set in ${\bf {\rm R}}^{n} $.  Although it is possible to apply  iterative-difference methods for solving a nonlinear problem (\ref{EQ__2_}), but it  is also possible to construct iterative methods that take into account the decomposition of the residual function. In this case, when solving nonlinear equations, methods (Shakhno et al. 2014, 2011; Shakhno 2016; C\u{a}tina\c{s} 1994; Iakymchuk et al. 2016) were constructed as combinations of the Newton method (Dennis et al. 1996; Ortega et al. 1970; Argyros 2008; Deuflhard 2004) and iterative-difference methods of chord (secant) and Kurchatov (Dennis et al. 1996; Ortega et al. 1970; Shakhno 2006, 2007; Argyros 2008; Ren et al. 2010, 2011; Shakhno et al. 2005).

In the paper (Shakhno 2017), we proposed a method for solving a nonlinear problem of least squares with a non-differentiable operator (\ref{EQ__2_}) constructed on the basis of the Gauss-Newton method method (Dennis et al. 1996; Ortega et al. 1970) and the Kurchatov type method (Shakhno et al. 2011, 2005; Ren 2011). We studied its local convergence under Lipschitz conditions and showed its effectiveness in comparison with other methods using test problems.

\section{Preliminaries}

To find the solution of the problem (\ref{EQ__2_}) we consider the Gauss-Newton-Kurchatov method (Shakhno 2017):
\begin{equation} \label{EQ__3_} 
\begin{array}{l} {x_{n+1} =x_{n} -(A_{n}^{{\bf \top }} A_{n} )^{-1} A_{n}^{{\bf \top }} (F(x_{n} )+G(x_{n} )),} \\ {} \\ {A_{n} =F'(x_{n} )+G(2x_{n} -x_{n-1} ,x_{n-1} ),{\rm \; \; \; \; \; \; }n=0,1,\ldots ,} \end{array} 
\end{equation} 
where $F'(x_{n} )$ is matrix of Jacobi of $F(x)$; $G(2x_{n} -x_{n-1} ,x_{n-1} )$ is the divided difference of the first order of functions (Ulm 1967), and  the points $2x_{n} -x_{n-1} ,x_{n-1} $; $x_{0} $, $x_{-1} $ are initial approximations. Method (\ref{EQ__3_}) is a combination of the  Gauss-Newton method  (Dennis et al. 1996; Ortega et al. 1970) and the Kurchatov type method  (Shakhno et al. 2011, 2005; Ren 2011).

If $m=n$, method (\ref{EQ__3_}) reduces to the Newton-Kurchatov method for solving the nonlinear equation $F(x)+G(x)=0$ (Shakhno et al. 2016, 2015; Hern\'{a}ndez-Ver\'{o}n 2017; Iakymchuk et al. 2016): 
\begin{equation} \label{EQ__4_} 
\begin{array}{l} {x_{n+1} =x_{n} -A_{n}^{-1} (F(x_{n} )+G(x_{n} )),} \\ {} \\ {A_{n} =F'(x_{n} )+G(2x_{n} -x_{n-1} ,x_{n-1} ),{\rm \; \; \; }n=0,1,\ldots } \end{array}.   
\end{equation} 
Setting in (\ref{EQ__3_}) $A_{n} =F'(x_{n} )+G(x_{n} ,x_{n-1} )$, we obtain a combination of the Gauss-Newton method (Dennis et al. 1996; Ortega et al. 1970) and the Secant type method (Ren et al. 2010; Shakhno et al. 2005) of the form (Shakhno et al. 2017) 
\begin{equation} \label{EQ__5_} 
\begin{array}{l} {x_{n+1} =x_{n} -A_{n}^{-1} (F(x_{n} )+G(x_{n} )),} \\ {} \\ {A_{n} =F'(x_{n} )+G(x_{n} ,x_{n-1} ),{\rm \; \; \; \; }n=0,1,\ldots } \end{array}.   
\end{equation} 
We need the following Lipschitz conditions.

\begin{definition}\label{d21} We say that the Fr\'echet derivative $F'$ satisfies the center Lipschitz conditions on $D$, if there exist  such that for each $x\in D$
\begin{equation} \label{EQ__6_} 
\left\|  F'(x)-F'(x^{*} )\right\| \le L_{0} \left\|  x-x^{*}  \right\| ,    
\end{equation} 
where $x^{*} \in D$ solves problem (\ref{EQ__2_}).
\end{definition}
 \begin{definition}\label{d22} We say that divided differences \textit{$G(\, \cdot \, ,\; \cdot )$}and\textit{ $G(\, \cdot \, ,\; \cdot \; ,\; \cdot \, )$ }satisfy the special Lipschitz conditions on $D\times D$ and $D\times D\times D$, if there exist $M_{0} >0$ and $N_{0} >0$ such that for each $x,y\in D$
\begin{equation} \label{EQ__7_} 
\left\|  G(x,y)-G(u,v) \right\| \le M_{0} {\rm (}\left\|  x-u \right\| +\left\|  y-v \right\| \, {\rm )},                                               
\end{equation} 
and
\begin{equation} \label{EQ__8_} 
\left\|  G(u,x,y)-G(v,x,y) \right\| \le N_{0} \left\|  u-v \right\| .  
\end{equation} 
\end{definition}
Let $B>0$ and $\alpha >0$. Define function  $h$ on $[0,\, +\infty )$ by
\begin{equation} \label{EQ__9_} 
h(t)\, =\, B\, {\bf [}(2\alpha +(L_{0} \, +2M_{0} )t+N_{0} t^{2} ]\, [(L_{0} /2+M_{0} )t+N_{0} t^{2} ]. 
\end{equation} 
Suppose that equation $h(t)=1$ has at least one positive solution. Denote by  $\gamma $ the smallest such solution. Set $D_{0} =D\bigcap \Omega (x^*,\gamma )$. 

\begin{definition}\label{d23} We say that the Fr\'echet  derivative $F'$ satisfies the restricted special Lipschitz conditions on $D_{0} $, if there exist  $L>0$ such that for euch $x,y\in D_{0} $
\begin{equation} \label{EQ__10_} 
\left\|  F'(x)-F'(y)\right\| \le L\left\|  x-y \right\|  
\end{equation} 
\end{definition}
\begin{definition}\label{d24} We say that divided differences \textit{$G(\, \cdot \, ,\; \cdot )$}and\textit{ $G(\, \cdot \, ,\; \cdot \; ,\; \cdot \, )$ }satisfy the special Lipschitz conditions on $D_{0} \times D_{0} $ and $D_{0} \times D_{0} \times D_{0} $, respectively, if there exist  $M>0$ and $N>0$ such that for each $x,y,u,v\in D_{0} $
\begin{equation} \label{EQ__11_} 
\left\|  G(x,y)-G(u,v) \right\| \le M{\rm (}\left\|  x-u \right\| +\left\|  y-v \right\| \, {\rm )} 
\end{equation} 
and
\begin{equation} \label{EQ__12_} 
\left\|  G(u,x,y)-G(v,x,y) \right\| \le N\left\|  u-v \right\| .   
\end{equation} 
\end{definition}
The following condition together with (\ref{EQ__7_}) and (\ref{EQ__8_}) have been used instead of the preceding ones in the study of such iterative methods (Shakhno 2017).

 \begin{definition}\label{d25} We say that the Fr\'echet derivative $F'$ satisfies the Lipschitz conditions on $D_{} $, if there exist  $L_{1} >0$ such that for euch $x,y\in D$
\begin{equation} \label{EQ__13_} 
\left\|  F'(x)-F'(y)\right\| \le L_{1} \left\|  x-y \right\|  
\end{equation} 
\end{definition}

Let $\Omega (x^{*} ,3r_{*} )={\rm \{ }x:\left\|  x-x^{*}  \right\| <3r_{*} {\rm \} }.$

\section{Convergence analysis  of the iterative process (\ref{EQ__3_})}
\setcounter{theorem}{0}

Next, we improve Theorem 1  (Shakhno et al. 2017).

 \begin{theorem} Let  function \textit{$F+G:{\bf {\rm R}}^{n} \to {\bf {\rm R}}^{m} $} be continuous on the open subset \textit{$D\subseteq {\bf {\rm R}}^{n} $},  \textit{$F$ }continuously differentiable in this domain, and let\textit{$G$ }be a continuous function. Assume that the problem (\ref{EQ__1_}) has a solution \textit{$x^{*} $ }in the domain and there exist the inverse operator $(A_{*}^{{\bf \top }} A_{*} )^{-1} ={\bf [}{\rm (}F'(x^{*} )+G(x^{*} ,x^{*} ){\rm )}^{{\rm \top }} {\rm (}F'(x^{*} )+G(x^{*} ,x^{*} ){\rm )}{\rm ]}^{-1} $
 and 
\[\left\|  (A_{*}^{{\bf \top }} A_{*} )^{-1}  \right\| \le B.\] 
Estimates (\ref{EQ__6_}), (\ref{EQ__7_}), (\ref{EQ__8_}), (\ref{EQ__10_}), (\ref{EQ__11_}), (\ref{EQ__12_}) hold and $\gamma $ given by (\ref{EQ__9_}) exists,
\begin{equation} \label{EQ__14_} 
\left\|  F(x^{*} )+G(x^{*} ) \right\| \le \eta ,{\rm \; \; \; \; \; \; }\left\|  F'(x^{*} )+G(x^{*} ,x^{*} ) \right\| \le \alpha ,  
\end{equation} 
\begin{equation} \label{EQ__15_} 
B(L+2M)\eta <1,  
\end{equation} 
\[\Omega (x^{*} ,3r_{*} )\subseteq D,       \] 
where $r_{*} $ is unique positive zero of the function $q$, given by 
 \begin{equation} \label{EQ__16_} \begin{array}{c} q(r) = B {\bf [}(\alpha +(L\, +2M)r+4Nr^{2} )((L/2+M)r+4Nr^{2} )+(L+2M+4Nr)\eta {\bf ]}\\
+B {\bf [}2\alpha +(L_{0} +2M_{0} )r+4N_{0} r^{2} {\bf ][}(L_{0} +2M_{0} )r+4N_{0} r^{2} {\bf ]}-1.   
\end{array}\end{equation} 

Then for $x_{0} ,\; x_{-1} \in \Omega (x^{*} ,r_{*} )$ the iterative process (\ref{EQ__3_}) is well defined, the sequence ${\rm \{ }x_{n} {\rm \} }$, $n=0,1,\ldots $, generated by it,  remains in the open subset $\Omega (x^{*} ,r_{*} )$, and converges to the solution $x^{*} $. Moreover,  the following error estimates hold for  $n=0,1,\ldots $ 
\begin{equation} \label{EQ__17_} \begin{array}{c}\left\|  x_{n+1} -x^{*}  \right\| \le C_{1} \left\|  x_{n} -x^{*}  \right\| +C_{2} \left\|  x_{n} -x_{n-1}  \right\| ^{2} +C_{3} \left\|  x_{n} -x^{*}  \right\| ^{2} \\
+C_{4} \left\|  x_{n-1} -x^{*}  \right\| ^{2} \left\|  x_{n} -x^{*}  \right\| ,  
\end{array}\end{equation} 
where
\[g(r)=B{\bf [}1-B{\rm (}2\alpha +(L_{0} +2M_{0} )r+4N_{0} r^{2} {\rm )(}(L_{0} +2M_{0} )r+4N_{0} r^{2} {\rm )}{\rm ]}^{-1} ,\] 
\[C_{1} =g(r_{*} )(L+2M)\eta ,{\rm \; \; \; \; \; }C_{2} =g(r_{*} )N\eta ,\] 
\begin{equation} \label{EQ__18_} 
C_{3} =g(r_{*} )(L{\bf /}2+M){\rm (}\alpha +(L+2M)r_{*} +4Nr_{*}^{2} {\rm )},  
\end{equation} 
\[C_{4} =g(r_{*} )N{\rm (}\alpha +(L+2M)r_{*} +4Nr_{*}^{2} {\rm )}. \] 
\end{theorem}

\textbf{Proof.} According to the intermediate value theorem on ${\bf [}0,\; r{\bf ]}$ the function $q$ for  a sufficiently large $r$ and by (\ref{EQ__15_}) has a positive zero denoted by $r_{*} $. But $q'(r)\ge 0$for $r\ge 0.$  So, this root is the only one on ${\bf [}0,\; r{\bf ]}$.

By assumption $x_{0} ,\; x_{-1} \in \Omega (x^{*} ,r_{*} )$. Then we have
\[\left\|  2x_{0} -\, x_{-1} -x^{*}  \right\| \le \left\|  x_{0} -x^{*}  \right\| +\left\|  x_{0} -\, x_{-1}  \right\|  \] 
\[\le \left\|  x_{0} -x^{*}  \right\| +\left\|   x_{0} -x^{*} \right\| +\left\|  x_{-1} -x^{*}  \right\| <3r_{*} .\] 
So, $2x_{0} -x_{-1} \in \Omega (x^{*} ,3r_{*} )$.

Let's denote $A_{n} =F'(x_{n} )+G(2x_{n} -x_{n-1} ,x_{n-1} )$. Let $n=0$ and we will get this estimate:
\begin{equation} \label{EQ__19_} \begin{array}{c}\left\|  I-(A_{*}^{{\bf \top }} A_{*} )^{-1} A_{0}^{{\bf \top }} A_{0}  \right\| =\left\|  (A_{*}^{{\bf \top }} A_{*}^{{\rm \; }} )^{-1} (A_{*}^{{\bf \top }} A_{*}^{{\rm \; }} -A_{0}^{{\bf \top }} A_{0} ) \right\| \\\\
=\left\|  (A_{*}^{{\bf \top }} A_{*}^{{\rm \; }} )^{-1} {\rm (}A_{*}^{{\rm \top }} (A_{*}^{{\rm \; }} -A_{0} )+(A_{*}^{{\rm \top }} -A_{0}^{{\rm \top }} )(A_{0} -A_{*}^{{\rm \; }} )\right. \\\\ 
\left. +(A_{*}^{{\bf \top }} -A_{0}^{{\bf \top }} )A_{*}^{{\rm \; }} {\rm )} \right\| \le \left\|  (A_{*}^{{\rm \top }} A_{*}^{{\rm \; }} )^{-1}  \right\| {\rm (}\, \left\|  A_{*}^{{\rm \top }}  \right\| \left\|  A_{*}^{{\rm \; }} -A_{0}  \right\| \\ \\
+\left\|  A_{*}^{{\bf \top }} -A_{0}^{{\bf \top }}  \right\| \left\|  A_{0} -A_{*}^{{\rm \; }}  \right\| +\left\|  A_{*}^{{\bf \top }} -A_{0}^{{\bf \top }}  \right\| \left\|  A_{*}^{{\rm \; }}  \right\| \, {\rm )}\\\\
\le B{\rm (}\alpha \left\|  A_{*}^{{\rm \; }} \, -\, A_{0}  \right\| +\left\|  A_{*}^{{\bf \top }} \, -\, A_{0}^{{\bf \top }}  \right\| \left\|  A_{0} \, -\, A_{*}^{{\rm \; }}  \right\| +\alpha \left\|  A_{*}^{{\bf \top }} \, -\, A_{0}^{{\bf \top }}  \right\| \, {\bf )}.  
\end{array}\end{equation} 
Using (\ref{EQ__8_}), we get
\begin{equation} \label{EQ__20_} \begin{array}{c}\left\|  G(2x_{0} -x_{-1} ,x_{-1} )-G(x_{0} ,x_{0} ) \right\| \\\\ 
=\left\|  G(2x_{0} -x_{-1} ,x_{-1} )-G(x_{0} ,x_{-1} )+G(x_{0} ,x_{-1} )-G(x_{0} ,x_{0} ) \right\| 
\\\\=\left\|  G(2x_{0} \, -x_{-1} ,x_{-1} ,x_{0} )(x_{0} \, -x_{-1} )-G(x_{0} ,x_{-1} ,x_{0} )(x_{0} \, -x_{-1} ) \right\| \\\\ 
\le N_{0} \left\|  x_{0} -x_{-1}  \right\| ^{2}  
\end{array}\end{equation} 
and
\begin{equation} \label{EQ__21_} \begin{array}{c}\left\|  G(2x_{0} -x_{-1} ,x_{-1} )-G(x_{0} ,x^{*} ) \right\| \\\\
=\left\|  G(2x_{0} -x_{-1} ,x_{-1} )-G(x_{0} ,x_{0} )+G(x_{0} ,x_{0} )-G(x_{0} ,x^{*} ) \right\| \\\\
\le N_{0} \left\|  x_{0} -x_{-1}  \right\| ^{2} +M_{0} \left\|  x_{0} -x^{*}  \right\| .  
\end{array}\end{equation} 
We use inequalities (\ref{EQ__7_}), (\ref{EQ__20_}), (\ref{EQ__21_}):
\begin{equation} \label{EQ__22_} \begin{array}{c}\left\|  A_{0} -A_{*}  \right\| =\left\|  (F'(x_{0} )+G(2x_{0} -x_{-1} ,x_{-1} ))-(F'(x^{*} )+G(x^{*} ,x^{*} )) \right\| \\\\
=\left\|  F'(x_{0} )-F'(x^{*} )+G(2x_{0} -x_{-1} ,x_{-1} )\right. \\\\
-\left. G(x_{0} ,x^{*} )+G(x_{0} ,x^{*} )-G(x^{*} ,x^{*} ) \right\| \\\\
\le L\left\|  x_{0} -x^{*}  \right\| +N\left\|  x_{0} -x_{-1}  \right\| ^{2} +2M\left\|  x_{0} -x^{*}  \right\| \\\\
=(L_{0} +2M_{0} )\left\|  x_{0} -x^{*}  \right\| +N_{0} \left\|  x_{0} -x_{-1}  \right\| ^{2} .  
\end{array}\end{equation}  
Then
\begin{equation} \label{EQ__23_} 
\left\|  A_{0}  \right\| \le \left\|  A_{0}  \right\| ||+\left\|  A_{0} -A_{*}  \right\| \le \alpha +(L_{0} +2M_{0} )\left\|  x_{0} -x^{*}  \right\| +N_{0} \left\|  x_{0} -x_{-1}  \right\| ^{2} .  
\end{equation} 
Then we obtain from the inequality (\ref{EQ__19_}) and the definition $r_{*} $ (\ref{EQ__16_})
\begin{equation} \label{EQ__24_} \begin{array}{c}\left\|  I-(A_{*}^{{\bf \top }} A_{*}^{{\rm \; }} )^{-1} A_{0}^{{\bf \top }} A_{0}  \right\| \le B\, {\bf [}2\alpha +(L_{0} +2M_{0} )\left\|  x_{0} -x^{*}  \right\| +N_{0} \left\|  x_{0} -x_{-1}  \right\| ^{2} {\bf ]} \\\\
\times {\bf [}(L_{0} +2M_{0} )\left\|  x_{0} -x^{*}  \right\| +N_{0} \left\|  x_{0} -x_{-1}  \right\| ^{2} {\bf ]}\\\\
\le B {\bf [}2\alpha +(L_{0} +2M_{0} )r_{*} +4N_{0} r_{*}^{2} {\bf ][}(L_{0} +2M_{0} )r_{*} +4N_{0} r_{*}^{2} {\bf ]}\\\\=h(r_{*} )<1. 
\end{array}\end{equation}  
By Banach's theorem on the inverse operator (Ortega et al. 270) there exists $(A_{0}^{{\bf \top }} A_{0} )^{-1} $  and we have from (\ref{EQ__24_})
\[\left\|  (A_{0}^{{\bf \top }} A_{0} )^{-1}  \right\| \le g_{0} =B{\rm \{ }1-B\, {\bf [}2\alpha +(L_{0} +2M_{0} )\left\|  x_{0} -x^{*}  \right\| +N_{0} \left\|  x_{0} -x_{-1}  \right\| ^{2} {\bf ]}\]  
\[\times {\bf [}(L_{0} +2M_{0} )\left\|  x_{0} -x^{*}  \right\| +N_{0} \left\|  x_{0} -x_{-1}  \right\| ^{2} {\bf ]}\, {\rm \} }^{-1} \]\[\le g(r_{*} )
=B{\rm \{ }1-B{\bf [}2\alpha +(L_{0} +2M_{0} )r_{*} +4N_{0} r_{*}^{2} {\bf ][}(L_{0} +2M_{0} )r_{*} +4N_{0} r_{*}^{2} {\bf ]}\, {\bf \} }^{-1} .\] 

Consequently,  iterate $x_{1} $ is well defined.

Then let's show that $x_{1} \in \Omega (x^{*} ,r_{*} )$. Using equality
\[A_{*}^{{\bf \top }} (F(x^{*} )+G(x^{*} ))=0 ,\] 

we will get an estimate
\[\left\|  x_{1} -x^{*}  \right\| =\] 
\[=\left\|  x_{0} -x^{*} -(A_{0}^{{\bf \top }} A_{0} )^{-1} (A_{0}^{{\bf \top }} (F(x_{0} )+G(x_{0} ))-A_{*}^{{\bf \top }} (F(x^{*} )+G(x^{*} )) \right\|  \] 
\[\le \left\|-(A_{0}^{{\bf \top }} A_{0} )^{-1}\right\| \left\| \, \left[\right. -A_{0}^{{\bf \top }} \left(\right.  A_{0} -\int _{0}^{1}F'(x^{*} +t (x_{0} -x^{*} )) dt\right. \] 
\[-\left. G(x_{0} ,x^{*} ))(x_{0} -x^{*} )+(A_{0}^{{\bf \top }} -A_{*}^{{\bf \top }} )(F(x^{*} )+G(x^{*} ))\left. \right] \, \right\| .\] 
Hence, taking into account (\ref{EQ__21_}), (\ref{EQ__23_}) and inequalities
\[\left\| \, A_{0} -\int _{0}^{1}F'(x^{*} +t (x_{0} -x^{*} ))dt-G(x_{0} ,x^{*} ) \right\| \] 
\[=\left\| \, F'(x_{0} )-\int _{0}^{1}F'(x^{*} +t (x_{0} -x^{*} ))\, dt+G(2x_{0} -x_{-1} ,x_{-1} )-G(x_{0} ,x^{*} ) \right\| \] 
\[=\left\| \,  \int _{0}^{1}(F'(x_{0} )-F'(x^{*} +t (x_{0} -x^{*} \, )))\, dt 
+G(2x_{0} -x_{-1} ,x_{-1} )-G(x_{0} ,x^{*} ) \right\|  \] 
\[\le \frac{1}{2} L\left\|  x_{0} -x^{*}  \right\| +M\left\|  x_{0} -x^{*}  \right\| +N\left\|   x_{0} -x_{-1}  \right\| ^{2}  \] 
\[\le \frac{1}{2} L\left\|  x_{0} -x^{*}  \right\| +M\left\|  x_{0} -x^{*}  \right\| +N{\rm (}\, \left\|  x_{0} -x^{*}  \right\| +\left\|  x_{-1} -x^{*}  \right\| \, {\rm )}^{2} \] 
we will get
\[\left\|  x_{1} -x^{*}  \right\| \le B\left\{\right. {\rm (}\alpha +(L+2M)\left\|  x_{0} -x^{*}  \right\| +N\left\|  x_{0} -x_{-1}  \right\| ^{2} {\rm )} \] 
\[\times \left(\right. \frac{1}{2} L\left\|  x_{0} -x^{*}  \right\| +M\left\|  x_{0} -x^{*}  \right\| +N\left\|  x_{0} -x_{-1}  \right\| ^{2} \left.\right)\left\|  x_{0} -x^{*}  \right\| \] 
\[+\eta \, {\rm (}(L+2M)\left\|  x_{0} -x^{*}  \right\| +N\left\|  x_{0} -x_{-1}  \right\| ^{2} \, {\rm )}\left. \right\} \] 
\[\times \left\{\right. 1-B\, {\bf [}2\alpha +(L+2M)\left\|  x_{0} -x^{*}  \right\| +N\left\|  x_{0} -x_{-1}  \right\| ^{2} {\bf ]}\] 
\[\times {\rm (}(L+2M)\left\|  x_{0} -x^{*}  \right\| +N\left\|  x_{0} -x_{-1}  \right\| ^{2} {\rm )}\left.  \right\}^{-1} \] 
\[=g_{0} \left\{\right.  {\rm (}\alpha +(L+2M)\left\|  x_{0} -x^{*}  \right\| +N\left\|  x_{0} -x_{-1}  \right\| ^{2} {\rm )} \] 
\[\times \left(\right. \frac{1}{2} L\left\|  x_{0} -x^{*}  \right\| +M\left\|  x_{0} -x^{*}  \right\| +N\left\|  x_{0} -x_{-1}  \right\| ^{2}  \left. \right)\left\|  x_{0} -x^{*}  \right\| \] 
\[+\eta \, {\rm (}(L+2M)\left\|  x_{0} -x^{*}  \right\| +N\left\|  x_{0} -x_{-1}  \right\| ^{2} {\rm )} \left. \right\}\] 
\[<g(r_{*} ){\bf [}{\rm (}\alpha +(L+2M)r_{*} +4Nr_{*}^{2} {\rm )(}(L/2+M)r_{*} +4Nr_{*}^{2} {\rm )}\] 
\[+(L+2M+4Nr_{*} )\eta {\bf ]}\, r_{*} =p(r_{*} )r_{*} =r_{*} ,\] 
where
\[p(r)=g(r)\, {\bf [}{\rm (}\alpha +(L+2M)r+4Nr^{2} {\rm )(}(L/2+M)r+4Nr^{2} {\rm )}\] 
\[+(L+2M+4Nr)\eta \, {\bf ]}.\] 
Hence, $x_{1} \in \Omega (x^{*} ,r_{*} )$ and inequality (\ref{EQ__16_}) is true for $n=0$.

 Assume that $x_{n} \in \Omega (x^{*} ,r_{*} )$ for $n=0,1,\ldots ,k$, and the estimate (\ref{EQ__17_}) for $n=0,1,\ldots ,k-1$, where $k\, \ge $1 is an integer, holds. Next we prove that $x_{n+1} \in \Omega (x^{*} ,r_{*} )$, and the estimate (\ref{EQ__17_}) holds for $n=k$.

 Define
\[\left\|  I-{\rm (}A_{*}^{{\bf \top }} A_{*}^{{\bf \top }} {\bf )}^{-1} A_{{}^{k} }^{{\bf \top }} A_{k}  \right\| =\left\|  {\bf (}A_{*}^{{\bf \top }} A_{*}^{{\rm \; }} {\bf )}^{-1} {\bf (}A_{*}^{{\bf \top }} A_{*}^{{\rm \; }} -A_{{}^{k} }^{{\bf \top }} A_{k} {\bf )} \right\| \] 
\[=\left\| \, {\rm (}A_{*}^{{\bf \top }} A_{*}^{{\rm \; }} {\bf )}^{-1} {\bf (}A_{*}^{{\bf \top }} (A_{*}^{{\rm \; }} -A_{k} )+(A_{*}^{{\bf \top }} -A_{{}^{k} }^{{\bf \top }} )(A_{k} -A_{*}^{{\rm \; }} )\right. \] 
\[+\left. {\rm (}A_{*}^{{\bf \top }} -A_{{}^{k} }^{{\bf \top }} {\bf )}A_{*} {\bf )}\, \right\| \le \left\|  {\bf (}A_{*}^{{\bf \top }} A_{*}^{{\rm \; }} {\bf )}^{-1}  \right\| \, {\bf (}\, \left\|  A_{*}^{{\bf \top }}  \right\| \left\|  A_{*}^{{\rm \; }} -A_{k}  \right\| \] 
\[+\left\|  A_{*}^{{\bf \top }} -A_{{}^{k} }^{{\bf \top }}  \right\| \left\|  A_{k} -A_{*}^{{\rm \; }}  \right\| +\left\|  A_{*}^{{\bf \top }} -A_{{}^{k} }^{{\bf \top }}  \right\| \left\|  A_{*}^{{\rm \; }}  \right\| \, {\rm )} \] 
\[\le B\, {\rm (}\alpha \left\|  A_{*}^{{\rm \; }} -A_{k}  \right\| +\left\|  A_{*}^{{\bf \top }} -A_{{}^{k} }^{{\bf \top }}  \right\| \left\|  A_{k} -A_{*}^{{\rm \; }}  \right\| +\alpha \left\|  A_{*}^{{\bf \top }} -A_{{}^{k} }^{{\bf \top }}  \right\| \, {\bf )} \] 
\[\le B\, {\bf [}2\alpha +(L+2M)\left\|  x_{k} -x^{*}  \right\| +N\left\|  x_{k} -x_{k-1}  \right\| ^{2} {\bf ]} \] 
\[\times {\bf [}(L/2+M)\left\|  x_{k} -x^{*}  \right\| +N\left\|  x_{k} -x_{k-1}  \right\| ^{2} {\bf ]} \] 
\[\le B\, {\bf [}2\alpha +(L+2M)r_{*} +4Nr_{*}^{2} {\bf ]}\, {\bf [}(L+2M)r_{*} +4Nr_{*}^{2} {\bf ]}=h(r_{*} )<1.\] 
So, ${\rm (}A_{{}^{k} }^{{\bf \top }} A_{k} {\bf )}^{-1} $ exists and
\[\left\|  (A_{k+1}^{{\bf \top }} A_{k+1} )^{-1}  \right\| \le g_{k} \, =B{\rm \{ }1-B\, {\bf [}2\alpha +(L_{0} \, +\, 2M_{0} )\left\|  x_{k} \, -\, x^{*}  \right\| +N_{0} \left\|  x_{k} \, -\, x_{k-1}  \right\| ^{2} {\bf ]} \] 
\[\times {\bf [}(L_{0} /2+M_{0} )\left\|  x_{k} -x^{*}  \right\| +N_{0} \left\|  x_{k} -x_{k-1}  \right\| ^{2} {\bf ]}{\rm \} }^{-1} \le g(r_{*} ).\] 

Therefore, the iteration $x_{k+1} $ is well defined, and we can get in turn $$\left\|  x_{k+1} -x^{*}  \right\| =\left\|  x_{k} -x^{*} -{\rm (}A_{k}^{{\bf \top }} A_{k} {\bf )}^{-1} {\bf (}A_{k}^{{\bf \top }} (F(x_{k} )+G(x_{k} ){\bf )}\right. $$
\[-\left. A_{*}^{{\bf \top }} (F(x^{*} )+G(x^{*} )) \right\| \le \left\|  -{\rm (}A_{k}^{{\rm \top }} A_{k} {\rm )}^{-1}  \right\|  \] 
\[\times \left\| \, \left[\right.  -A_{k}^{{\bf \top }} \left(\right. A_{k} -\int _{0}^{1}F'(x^{*} +t (x_{k} -x^{*} ))\, dt\right. \] 
\[-\left. G(x_{k} ,x^{*} ) \left. \right)(x_{k} -x^{*} )+{\rm (}A_{k}^{{\bf \top }} -A_{*}^{{\bf \top }} {\bf )}(F(x^{*} )+G(x^{*} ))\left. \right]^{} \right\|  \] 
\[\le g_{k} {\rm \{ }\, {\bf [}\alpha +(L+2M)\left\|  x_{k} -x^{*}  \right\| +N\left\|  x_{k} -x_{k-1}  \right\| ^{2} {\bf ]} \] 
\[\times {\bf [}(L/2+M)\left\|  x_{k} -x^{*}  \right\| +N\left\|  x_{k} -x_{k-1}  \right\| ^{2} {\bf ]}\left\|  x_{k} -x^{*}  \right\| \] 
\[+\eta \, {\rm (}(L+2M)\left\|  x_{k} -x^{*}  \right\| +N\left\|  x_{k} -x_{k-1}  \right\| ^{2} {\rm )}{\rm \} } \] 
\[\le g(r_{*} ){\rm \{ }\, {\bf [}\alpha +(L+2M)\left\|  x_{k} -x^{*}  \right\| +N\left\|  x_{k} -x_{k-1}  \right\| ^{2} {\bf ]} \] 
\[\times {\bf [}(L/2+M)\left\|  x_{k} -x^{*}  \right\| +N\left\|  x_{k} -x_{k-1}  \right\| ^{2} {\bf ]}\left\|  x_{k} -x^{*}  \right\| \] 
\[+\eta \, {\rm (}(L+2M)\left\|  x_{k} -x^{*}  \right\| +N\left\|  x_{k} -x_{k-1}  \right\| ^{2} {\rm )\} }<p(r_{*} )r_{*} =r_{*} ,\] 
i.e. $x_{k+1} \in \Omega (x^{*} ,r_{*} )$, and estimate (\ref{EQ__17_}) holds for $n=k$

Consequently, the iterative process (\ref{EQ__3_}) is well defined, $x_{n} \in \Omega (x^{*} ,r_{*} )$ for all $n\ge 0$, and estimate (\ref{EQ__17_}) holds for all $n\ge 0$.

Next, we prove that $x_{n} \to x^{*} $ for $n\to \infty $. Define  functions $a$ and $b$ on ${\bf [}0,r_{*} {\bf ]}$ by:
\begin{equation} \label{EQ__25_}\begin{array}{c} a(r)=g(r)((L+2M+3Nr)\eta +\varphi (r)((L/2+M)r+4Nr^{2} )),\\\\ 

b(r)=g(r)Nr \eta ,  
\end{array}\end{equation} 
where $\varphi (r)=\alpha +(L+2M)r+4Nr^{2} $.

According to the choice $r_{*} $, we have 
\begin{equation} \label{EQ__26_} 
a(r_{*} )\ge 0,{\rm \; \; \; \; }b(r_{*} )\ge 0,{\rm \; \; \; \; }a(r_{*} )+b(r_{*} )=1.  
\end{equation} 
Using the estimate (\ref{EQ__17_}), the definition of constants $C_{i} ,\, \, i=1,2,3,4$, as well as the functions $a$ and $b$, for $n\ge 0$, we obtain
\begin{equation} \label{EQ__27_}\begin{array}{c} \left\|  x_{n+1} -x^{*}  \right\| \le (C_{1} +C_{3} r+4C_{4} r_{*}^{2} )\left\|  x_{n} -x^{*}  \right\| +C_{2} {\rm (}\, \left\|  x_{n} -x^{*}  \right\| ^{2} \\\\ 
+2\left\|  x_{n-1} -x^{*}  \right\| \left\|  x_{n} -x^{*}  \right\| +\left\|  x_{n-1} -x^{*}  \right\| ^{2} {\rm )}\\\\<(C_{1} +3C_{2} r_{*} 
+C_{3} r_{*} +4C_{4} r_{*}^{2} )\left\|  x_{n} -x^{*}  \right\| +C_{2} r_{*} \left\|  x_{n-1} -x^{*}  \right\| \\\\ 
=a(r_{*} )\left\|  x_{n} -x^{*}  \right\| +b(r_{*} )\left\|  x_{n-1} -x^{*}  \right\| .  
\end{array}\end{equation} 

Similarly to (Ren at al. 2011), we prove that under the conditions (\ref{EQ__25_}), (\ref{EQ__26_}) the sequence ${\rm \{ }x_{n} {\rm \} }$for $n\to \infty $converges to $x^{*} $.

First of all, for a real number $r_{*} >0$ and initial points $x_{0} ,x_{-1} \in \Omega (x^{*} ,r_{*} )$ there exists a real number $r'$ such that $0<r'<r_{*} $, $x_{0} ,x_{-1} \in \Omega (x^{*} ,r')$. Then all the above estimates for the sequence ${\rm \{ }x_{n} {\rm \} }$ are valid, if replaced $r_{*} $ by $r'$. In particular, from (\ref{EQ__27_})  for $n\ge 0$, we get
\begin{equation} \label{EQ__28_} 
\left\|  x_{n+1} -x^{*}  \right\| \le a\left\|  x_{n} -x^{*}  \right\| +b\left\|  x_{n-1} -x^{*}  \right\| ,  
\end{equation} 
where $a=a(r')$, $b=b(r')$.

Clearly, we  also have
\[a\ge 0,{\rm \; \; \; \; }b\ge 0,{\rm \; \; \; \; }a+b<a(r_{*} )\left\|  x_{n} -x^{*}  \right\| +b(r_{*} )\left\|  x_{n-1} -x^{*}  \right\| <1.\] 

Define sequences ${\rm \{ }\theta _{n} {\rm \} }$, ${\rm \{ }\rho _{n} {\rm \} }$:
\begin{equation} \label{EQ__29_}\begin{array}{c} \theta _{n} =\displaystyle\frac{\left\|  x_{n} -x^{*}  \right\| }{r'} ,{\rm \; \; \; \; \; \; }n=-1,0,1,\ldots ,\\\\
 
\rho _{-1} =\theta _{-1} ,{\rm \; \; \; }\rho _{0} =\theta _{0} ,{\rm \; \; \; }\rho _{n+1} =a\rho _{n} +b\rho _{n} {}_{-1} ,{\rm \; \; \; }n=0,1,2,\ldots .  
\end{array}\end{equation} 
We divide the two parts of inequality (\ref{EQ__28_}) into $r'$  and obtain $\theta _{n+1} =a\theta _{n} +b\theta _{n} {}_{-1} ,{\rm \; \; \; \; }n=0,1,2,\ldots $.

 By definition of the sequence ${\rm \{ }\rho _{n} {\rm \} }$, we have 
\begin{equation} \label{EQ__30_} 
0\le \theta _{n} \le \rho _{n} ,{\rm \; \; \; \; }n=-1,0,1,\ldots .  
\end{equation} 
For the sequence ${\rm \{ }\rho _{n} {\rm \} }$ known explicit formulas 
\begin{equation} \label{EQ__31_} 
\rho _{n} =\omega _{1} \lambda _{1}^{n} +\omega _{2} \lambda _{2}^{n} ,{\rm \; \; \; \; }n=-1,0,1,\ldots ,                     
\end{equation} 
where
\[\lambda _{1} =\frac{a-\sqrt{a^{2} +4b} }{2} ,{\rm \; \; \; \; \; }\lambda _{2} =\frac{a+\sqrt{a^{2} +4b} }{2} \] 
and
\[\omega _{1} =\frac{\lambda _{2}^{-1} \rho _{0} -\rho _{-1} }{\lambda _{2}^{-1} -\lambda _{1}^{-1} } ,{\rm \; \; \; \; \; }\omega _{2} =\frac{\rho _{-1} -\lambda _{1}^{-1} \rho _{0} }{\lambda _{2}^{-1} -\lambda _{1}^{-1} } .\] 
Note that
\[0\le \left| \lambda _{1}  \right|\le \left| \lambda _{2}  \right|<\frac{a+\sqrt{a^{2} +4(1-a)} }{2} =\frac{a+2-a}{2} =1.\] 
Taking into account (\ref{EQ__30_}) and (\ref{EQ__31_}), we conclude that ${\rm \{ }\theta _{n} {\rm \} }\to 0$ as $n\to \infty $. Therefore,  we conclude that $x_{n} \to x^{*} $ as $n\to \infty $.     
\hfill ${\rlap{$\sqcup$}\displaystyle\sqcap} $
 \begin{remark} If $L_{0} =L=L_{1} $, $M_{0} =M$ and $N_{0} =N$, our results specialize  to the corresponding ones  (Shakhno 2017). Otherwise they constitute an improvement. As an example let $q_{1} ,\, g_{1} ,\, C_{1}^{1} ,\, C_{2}^{1} ,\, C_{3}^{1} ,\, C_{4}^{1} ,\, r_{*}^{1} $ used in (Shakhno 2017) denote the functions and parameters, where $L_{0} ,\, L,\, M,\, N$ are replaced by $L_{1} ,\, L_{1} ,\, M_{0} ,\, N_{0} $, respectively. Then, since $L_{0} \le \, L_{1} $, $L\le \, L_{1} \, $, $M\le \, M_{0} \, $, $N\le \, N_{0} \, $and since $D_{0} \, \subseteq D$, we have $q(r)\le q_{1} (r)$,  $g(r)\le g_{1} (r)$,  $C_{1} \le C_{1}^{1} $,  $C_{2} \le C_{2}^{1} $,  $C_{3} \le C_{3}^{1} $, $C_{4} \le C_{4}^{1} $, so $r_{*}^{1} \le r_{*} $,  and the new error bounds are tighter than the corresponding ones (\ref{EQ__23_}) (Shakhno 2017) . 

 Moreover, we have
\[B(L_{0} +2M_{0} )\eta <1\, \, \, \, \, \, \, \, \Rightarrow \, \, \, \, \, \, \, \, \, B(L+2M)\eta <1\, \] 
but not vice versa, unless if $L_{0} =L$ and $M_{0} =M$. 

 Hence, the new sufficient convergence criteria for method (\ref{EQ__3_}) are weaker. These advantages are obtained under the same computational cost as (Shakhno 2017), since in practice the new constants are special cases of the previous ones.
\end{remark}
\begin{corollary} In the case of zero residual, the convergence order of the iterative process (\ref{EQ__3_}) is quadratic.
\end{corollary}
 If $\eta =0$, we have a nonlinear least squares problem with zero residual in the solution. Then the constants $C_{1} =0$ and $C_{2} \, =0$ and (\ref{EQ__17_}) reduces to 
\begin{equation} \label{EQ__32_} 
\left\|  x_{n+1} -x^{*}  \right\| \le C_{3} \left\|  x_{n} -x^{*}  \right\| ^{2} +C_{4} \left\|  x_{n} -x_{n-1}  \right\| ^{2} \left\|  x_{n} -x^{*}  \right\| .             
\end{equation} 
It follows from the inequality (\ref{EQ__32_}) that the order of convergence (\ref{EQ__3_}) is not higher than quadratic. Consequently, there exist a constant $C_{5} \ge 0$ and a positive integer $N$ such that for all $n\ge N$
\[\left\|  x_{n} -x^{*}  \right\| \ge C_{5} \left\|  x_{n-1} -x^{*}  \right\| ^{2} .\] 
By
\[\left\|  x_{n} -x^{*}  \right\| \le \left\|  x_{n-1} -x^{*}  \right\| ,\] 
we have
\[\left\|  x_{n} -x_{n-1}  \right\| ^{2} \le {\rm (} \left\|  x_{n} -x^{*}  \right\| +\left\|  x_{n-1} -x^{*}  \right\| \, {\rm )}^{2} \le 4\left\|  x_{n-1} -x^{*}  \right\| ^{2} ,\] 
and from (\ref{EQ__32_}) we have
\[\left\|  x_{n+1} -x^{*}  \right\| \le C_{3} \left\|  x_{n} -x^{*}  \right\| ^{2} +4C_{4} \left\|  x_{n-1} -x^{*}  \right\| ^{2} \left\|  x_{n} -x^{*}  \right\|  \] 
\[\le C_{3} \left\|  x_{n} -x^{*}  \right\| ^{2} +4\frac{C_{4} }{C_{5} } \left\|  x_{n} -x^{*}  \right\| ^{2} =C_{6} \left\|  x_{n} -x^{*}  \right\| ^{2} .\] 

Consequently, the convergence order  of the iterative process (\ref{EQ__3_}) is quadratic.

As we see from the estimates (\ref{EQ__17_}) and (\ref{EQ__18_}), the convergence of the iterative process (\ref{EQ__3_}) essentially depends on the terms containing the values $\eta $, $\alpha $, $L$, $M$ and $N$.

 For problems with zero residual in the solution ($\eta =0$), the quadratic convergence of the iterative process (\ref{EQ__3_}) is established.

For problems with a small residual in the solution ($\eta $ -- "small") and with weak nonlinearity ($\alpha $, $L_{0} $, $L$, $M$ and  $N$-- "small"), the convergence of the iterative process is linear. In the case of large residual ($\eta $ -- "large") or for strongly nonlinear problems ($\alpha $,  $L_{0}$, $L$, $M$ and $N$ -- "large"), the iterative process (\ref{EQ__3_}) may not converge at all.

\section{Results of numerical experiment}

On several test cases, we compare the convergence rates of the Gauss-Newton-Kurchatov method (\ref{EQ__3_}), the Gauss-Newton-Secant method (\ref{EQ__5_} ) and the Secant-type difference method\textbf{ }(Ren et al. 2010; Shakhno et al 2005) \textbf{}
\begin{equation} \label{EQ__33_}\begin{array}{c}  x_{n+1} =x_{n} -{\rm (}A_{n}^{{\bf \top }} A_{n} {\bf )}^{-1} A_{n}^{{\bf \top }} (F(x_{n} )+G(x_{n} )),\\\\
A_{n} =F(x_{n} ,x_{n-1} )+G(x_{n} ,x_{n-1} ),{\rm \; \; \; \; \; }n=0,1,\ldots ,  
\end{array}\end{equation} 
and the Kurchatov-type difference method\textbf{ }(Ren et al. 2011;\textbf{ }Shakhno et al. 2005) 
\begin{equation} \label{EQ__34_}\begin{array}{c}  x_{n+1} =x_{n} -{\rm (}A_{n}^{{\bf \top }} A_{n} {\bf )}^{-1} A_{n}^{{\bf \top }} (F(x_{n} )+G(x_{n} )),\\\\
A_{n} =F(2x_{n} -x_{n-1} ,x_{n-1} )+G(2x_{n} -x_{n-1} ,x_{n-1} ),{\rm \; \; \; \; }n=0,1,\ldots .  
\end{array}\end{equation} 
We tested methods on nonlinear systems with a non-differentiable operator with zero and non-zero residuen. The classical Gauss-Newton method and the Newton method cannot apply to  solving these problems.

 Solution results are with accurate $\varepsilon =10^{-8} $. The additional approximation was chosen as follows: $x_{-1} =x_{0} -10^{-4} $. The calculations were carried out until the conditions were fulfilled 

 $$\left\|  x_{n+1} -x_{n}  \right\| \le \varepsilon  \quad \rm{and} \quad  \left\|  A_{n}^{{\bf \top }} (F(x_{n} )+G(x_{n} )) \right\| \le \varepsilon ,$$
 with $f(x)={\mathop{\min }\limits_{x\in {\bf {\rm R}}^{n} }} \frac{1}{2} (F(x)+G(x))^{{\bf \top }} (F(x)+G(x))$.

\vspace{2mm}
\textbf{Example} \textbf{\textit{1}}\textit{  }(Shakhno et al. 2014; Argyros 2008; C\u{a}tina\c{s} 1994):
\[\left\{\begin{array}{l} {3x^{2} y+y^{2} -1+\left| x-1 \right|=0,} \\ {x^{4} +xy^{3} -1+\left| y \right|=0,} \end{array}\right. \] 
\[(x^{*} ,y^{*} )\approx (0.89465537,{\rm \; }0.32782652), f(x^{*} )=0.\] 

\textbf{Example\textit{ 2.} }$n=2,{\rm \; \; }m=3$:
\[\left\{\begin{array}{l} {3x^{2} y+y^{2} -1+\left| x-1 \right|=0,} \\ {x^{4} +xy^{3} -1+\left| y \right|=0,} \\ {\left| x^{2} -y \right|={\rm 0,}} \end{array}\right. \] 
\[(x^{*} ,y^{*} )\approx (0.74862800,{\rm \; }0.43039151), f(x^{*} )\approx 4.0469349\cdot 10^{-2} .\] 

Table 1 shows the results of a numerical experiment. In particular, the investigated methods are compared by the number of iterations performed to find a solution with a given accuracy.

\vspace{2mm}
{\it Table 1}. Number of iterations for solving of the test problems

\begin{tabular}{|p{0.5in}|p{0.6in}|p{0.6in}|p{0.8in}|p{0.5in}|p{0.7in}|} \hline 
Example  &  $(x_{0} ,y_{0} )$ & Kurchatov type's method (\ref{EQ__34_}) & Gauss-Newton-Kurchatov method (\ref{EQ__3_}) & Secant type's method (\ref{EQ__33_}) & Gauss-Newton-Secant method (\ref{EQ__5_}) \\ \hline 
\textbf{\textit{1}} & (1, 0.1) & 6 & 5 & 6 & 5 \\ 
\cline {2-6}
 & (3, 1) & 12 & 9 & 11 & 10 \\ \cline {2-6}
 & (0.5, 0.5) & 12 & 10 & 18 & 10 \\ \hline 
\textbf{\textit{2}} & (1, 0.1) & 16 & 14 & 21 & 11 \\\cline {2-6} 
 & (3, 1) & 21 & 18 & 25 & 15 \\ \cline {2-6}
 & (0.5, 0.5) & 16 & 14 & 19 & 13 \\\hline 
\end{tabular}

\section{Conclusions}

It follows from the theoretical results, practical calculations and comparison of the results obtained, that the combined differential-difference methods (\ref{EQ__3_}) and (\ref{EQ__5_}) converge faster than the Kurchatov type method (\ref{EQ__34_}) and the Secant type method  (\ref{EQ__33_}). As proved, in the case of zero residual, method (\ref{EQ__3_}) has a quadratic order of convergence and does not require the calculation of derivatives from a non-differentiable part of the operator. Then the method (\ref{EQ__3_}), as well as the method (\ref{EQ__5_}), are effective methods for solving   nonlinear  least squares problems with  non-differentiable operator.

\end{document}